\input amstex
\documentstyle{amsppt}
\input epsf
\magnification 1200
\vcorrection{-10mm}
\NoBlackBoxes

\def\letterwise{\equiv}
\def\Z{\Bbb Z}
\def\B{\bold B}

\rightheadtext{On the Hurwitz action on quasipositive 3-braids }
\topmatter
\title
            On the Hurwitz action on quasipositive factorizations of 3-braids
\endtitle
\author
            S.~Yu.~Orevkov
\endauthor
\address   Steklov Mathematical Institute, Moscow
\endaddress
\address   Universit\'e Paul Sabatier (Toulouse-3)
\endaddress
\email orevkov\@math.ups-tlse.fr
\endemail
\endtopmatter

\vskip -5mm
\document
Let $\B_3$ be the group of braids with three strings:
$\B_3=\langle\sigma_1,\sigma_2\mid\sigma_1\sigma_2\sigma_1=
\sigma_2\sigma_1\sigma_2\rangle$. A
{\it quasipositive factorization} of $X\in\B_3$
is a collection $(X_1,\dots,X_k)\in\B_3^k$
such that $X_1\dots X_k=X$, and each of the
$X_i$'s is conjugate in $\B_3$ to the standard generator $\sigma_1$.
Note that $\sigma_1$ and $\sigma_2$ are conjugate to each other in $\B_3$.
A braid $X$ is called {\it quasipositive} if it admits at least one
quasipositive factorization.

Let $G$ be a group. We define the mappings $\Sigma_i:G^k\to G^k$,
$1\le i<k$, by setting $\Sigma_i(X_1,\dots,X_k)=(Y_1,\dots,Y_k)$ where
$(Y_i,Y_{i+1})=(X_iX_{i+1}X_i^{-1},X_i)$ and $Y_j=X_j$ for $j\not\in\{i,i+1\}$.
These mappings are invertible because $(X_i,X_{i+1})=(Y_{i+1},Y_{i+1}^{-1}Y_iY_{i+1})$.
If $(X_1,\dots,X_k)$ is a quasipositive factorization of a braid $X$,
then it is easy to see that $\Sigma_i(X_1,\dots,X_k)$ is also
a quasipositive factorization of the same braid.
The correspondence $\sigma_i\mapsto\Sigma_i$ is an action of the braid group
$\B_k$ on the set $G^k$. This action is called the
{\it Hurwitz action}. Factorizations belonging to the same orbit of the
Hurwitz action are called {\it Hurwitz equivalent}.

If $U$, $V$ are words over the alphabet $\Cal A=\{\sigma_1,\sigma_2\}$ (such words,
as well as the braids represented by them, are called {\it positive\/}), then
$U\letterwise V$ stands for letterwise coincidence and $U=V$ stands for equality
in $\B_3$. We set $\Delta=\sigma_1\sigma_2\sigma_1$.

If $W\letterwise a_1\dots a_n$ is a positive word and $I=\{i_1,\dots,i_k\}$,
$1\le i_1<\dots<i_k\le n$, then we set
$W\setminus I \letterwise W_1\dots W_{k-1}$ and
$W_I = (A_1 x_1 A_1^{-1},\,\dots,\,A_k x_k A_k^{-1})$
where $x_j=a_{i_j}$, $W_1,\dots,W_{k+1}$ are pieces into which
$W$ is split by $x_1,\dots,x_k$
(i.~e., $W_j \letterwise a_l\dots a_m$ where $l=i_{j-1}+1$, $m=i_j-1$, $i_0=0$, $i_{k+1}=n+1$)
and $A_j=W_1\dots W_j$.
It is easy to see that $W_I$ is a quasipositive factorization of $W\Delta^{-p}$
if and only if $W\setminus I = \Delta^p$.
If $i-1\not\in I$, $i\in I$, and $a_{i-1}=a_i$, then
$W_I=W_{\{i-1\}\cup(I\setminus\{i\})}$. When $I$ does not contain such $i$, we say that
$I$ is {\it $W$-minimal}.

\proclaim{ Theorem 1 }
Let $W$ be a positive word and $p\ge0$. Suppose that the braid
$X=W \Delta^{-p}$ is quasipositive.
Then every orbit of the Hurwitz action on quasipositive factorizations of $X$
contains an element of the form $W_I$ with $W$-minimal $I$.
\endproclaim

\proclaim{ Corollary 1 } {\rm(see [3]).} A braid $X\in\B_3$
is quasipositive if and only if any positive word $W$ such that
$X=W\Delta^{-p}$ admits a removal of some letters so that the remaining word
is equal to $\Delta^p$ in $\B_3$.
\qed\endproclaim

Since any 3-braid can be presented as $W\Delta^{-p}$, Theorem 1
provides an algorithm for finding representatives of all orbits of the
Hurwitz action. The algorithm can be optimized in the `branch and bound' spirit
by analogy with [4; \S6].

\proclaim{ Corollary 2 }
The number of orbits of the Hurwitz action for any 3-braid is finite.
\footnote{(Added in 2019) There is a mistake here in the published version of
this article.}
\qed\endproclaim

\proclaim{ Corollary 3 }
Any two quasipositive factorizations of a
{\rm positive} 3-braid are Hurwitz equivalent.
\endproclaim

\demo{ Proof }
\footnote{(Added in 2019) In fact it is enough to observe that $I$ is unique.}
By Theorem 1, it is enough to check that the factorizations
$(\sigma_1,\sigma_2,\sigma_1)$ and $(\sigma_2,\sigma_1,\sigma_2)$
are Hurwitz equivalent. Indeed,
$(\sigma_1,\,\sigma_2,\,\sigma_1)
   \overset{\Sigma_2}\to\mapsto
   (\sigma_1,\,\sigma_2\sigma_1\sigma_2^{-1},\,\sigma_2)
   \overset{\Sigma_1}\to\mapsto
$\par\noindent$
   (\sigma_1\sigma_2\sigma_1\sigma_2^{-1}\sigma_1^{-1},\,\sigma_1,\,\sigma_2)
   = (\sigma_2,\,\sigma_1,\,\sigma_2)$.
\qed
\enddemo

The Birman--Ko--Lee (BKL) presentation for $\B_3$ is
$$
  \B_3=\langle\sigma_0,\sigma_1,\sigma_2\mid
  \sigma_2\sigma_1=\sigma_1\sigma_0=\sigma_0\sigma_2\rangle                  \eqno(2)
$$
where $\sigma_1$ and $\sigma_2$ are the same as above, and hence
$\sigma_0=\sigma_1^{-1}\sigma_2\sigma_1$. Words over the alphabet $\{\sigma_0,\sigma_1,\sigma_2\}$
and the braids represented by them will be called {\it BKL-positive}.
Any 3-braid can be written in the form $W\delta^{-p}$ with a BKL-positive $W$ and
$\delta=\sigma_2\sigma_1$.

\proclaim{ Theorem 2 }
Let $W$ be a BKL-positive word and $p\ge0$. Suppose that the braid
$X=W \delta^{-p}$ is quasipositive.
Then every orbit of the Hurwitz action on quasipositive factorizations of $X$
contains an element of the form $W_I$ with $W$-minimal $I$.
\endproclaim

\proclaim{ Corollary 4 }
Any two quasipositive factorizations of a
{\rm BKL-positive} 3-braid are Hurwitz equivalent.
\endproclaim

Of course, a BKL analog of Corollary 1 holds as well.
In spite of the similarity between Theorems 1 and 2, their proofs are very different.
Our proof of Theorem 1 is more geometric. It is inspired by the proof of the main
result in [2]. The proof of Theorem 2 is purely combinatorial, it is
in the spirit of [3].

Let $e:\B_3\to\Z$ be the group homomorphism such that $e(\sigma_1)=e(\sigma_2)=1$.

\proclaim{ Theorem 3 } If $X\in\B_3$ and
$e(X)=2$, then $X$ has at most two orbits of the Hurwitz action.
\endproclaim

\definition{\bf Example } Let $W\letterwise\sigma_1^2\sigma_2^2\sigma_1^2\sigma_2^2$. Then
$W_{\{1,5\}}\not\sim W_{\{3,7\}}$. Hence, by Theorems 1 and 3, the braid $W\Delta^{-2}$ has
exactly two orbits of the Hurwitz action.
\enddefinition

\smallskip\noindent
{\bf Remark.} Theorem 2 and its proof extend without changes to the case of
Artin-Tits groups of type $I_2 (p)$ if one defines BKL-positive words as positive
words in the generators of the presentation $\langle a_1,\dots, a_p \mid
 a_p a_{p-1} = a_{p-1} a_{p-2} = \dots = a_2 a_1 = a_1 a_p\rangle$.


\subhead \S1. Admissible graphs and quasipositive factorizations
\endsubhead
We fix a disk $D$ and a point $q$ on its boundary $\partial D$.
Let $\Gamma$ be an oriented graph embedded in $D\setminus\{q\}$
and let the numbers 1 and 2 be assigned to every edge of $\Gamma$.
Let $V(\Gamma)$ be the set of all vertices of
$\Gamma$ and let $V_n(\Gamma)$ be the set of the vertices adjacent to $n$ edges.
We set $\partial\Gamma=\Gamma\cap\partial D$, $R(\Gamma)=V_6(\Gamma)$,
and $B(\Gamma)=V_1(\Gamma)\setminus\partial\Gamma$. Elements of $B(\Gamma)$
will be called {\it branch points}.

\midinsert
\centerline{
  \epsfxsize=15mm\epsfbox{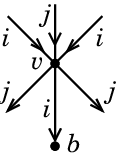}   \hskip 20mm
  \epsfxsize=25mm\epsfbox{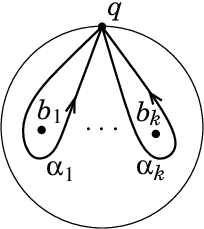}
}
\smallskip
\centerline{
  {\bf Fig.~1.} $i\ne j$ \hskip23mm
  {\bf Fig.~2.} \hskip10mm
}
\endinsert

We say that $\Gamma$ is an {\it admissible graph\/}
if $\partial\Gamma\subset V_1(\Gamma)$, $V(\Gamma) = V_1(\Gamma)\cup R(\Gamma)$,
and the edges incident to any vertex $v\in R(\Gamma)$ are oriented and labeled
as in Fig.~1. An admissible graph $\Gamma$ will be called {\it quasipositive\/}
if the edge incident to any branch point $b$ is oriented towards $b$ (see Fig.~1).

If $\alpha$ is a path in $D\setminus V(\Gamma)$ transversal to the edges
of $\Gamma$, then we define the word $\Gamma(\alpha)$ as
$\sigma_{i_1}^{\pm1}\sigma_{i_2}^{\pm1}\dots$ where $i_1,i_2,\dots$ are the
numbers assigned to the edges successively crossed by $\alpha$ and the signs are
chosen according to the orientation of these edges, so that a positively oriented loop
around a branch point labeled by $i$ corresponds to $\sigma_i$.
The word corresponding to the positive circuit
along $\partial D$ starting from $q$ will be called the {\it boundary word\/} of
$\Gamma$. We denote it by $\Gamma(\partial D)$.
If $\Gamma$ is an admissible graph, then it is easy to see that paths which are
homotopic in $D\setminus B(\Gamma)$ define the same braid.

Every quasipositive graph $\Gamma$ uniquely determines a Hurwitz equivalence class
of quasipositive factorizations of the boundary braid as follows. Let
$B(\Gamma)=\{b_1,\dots,b_k\}$. We choose pairwise distinct paths $\alpha_1,\dots,\alpha_k$
as in Fig.~2. Then $(\Gamma(\alpha_1),\,\dots,\,\Gamma(\alpha_k))$ is a
quasipositive factorization of $\Gamma(\partial D)$. The collection of paths $(\alpha_i)$
is defined up to a diffeomorphism of the disk identical on the boundary, i.~e., up to
the action of the braid group $\B_k$. One can check
that this is the Hurwitz action.

\proclaim{ Lemma 1 } Any quasipositive factorization of a given word can be
represented by a quasipositive graph.
\endproclaim

\demo{ Proof }
If words $W_1$ and $W_2$ are equal in $\B_3$ and
a disk $D_2$ is inside $D_1$, then there exists a branch point free admissible graph
$\Gamma$ in the annulus $D_1\setminus D_2$ such that
$\Gamma(\partial D_i)=W_i$, $i=1,2$.
It is enough to check this fact when $W_2$ is obtained from $W_1$ either by a braid group
relation or by inserting or removing
$\sigma_i^{\pm1}\sigma_i^{\mp1}$. Let us construct a graph which defines
a given quasipositive factorization $(X_1,\dots,X_k)$
of a given word $W$, $X_i=a_i\sigma_1 a_i^{-1}$. We consider nested disks
$D_3\subset D_2\subset D_1\subset D$. In
$D\setminus D_1$, we construct a graph realizing the equality $W=X_1\dots X_k$.
In $D_1\setminus D_2$, we complete the edges corresponding to the central
$\sigma_1$'s  by adding branch points. Finally, in $D_2\setminus D_3$, we
realize the equality $(a_1a_1^{-1})\dots(a_ka_k^{-1})=1$.
\qed\enddemo

\proclaim{ Lemma 2 }
If two quasipositive graphs coincide outside a disk
$U\subset D$ and if each of them has at most one branch point in $U$, then the graphs
define the same quasipositive factorization of the boundary braid.
\qed\endproclaim

\midinsert
\centerline{ \epsfxsize=75mm\epsfbox{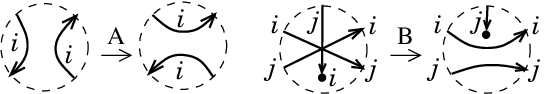} }
\smallskip
\centerline{ \bf Fig.~3. $i\ne j$ }
\endinsert

In fact, we need Lemma 2 only in two cases: when one graph is obtained from the other
by the modifications shown in Fig.~3.


\subhead \S2. Proof of Theorem 1
\endsubhead
Let $\Delta_{-p}$ be the word $\sigma_1^{-1}\sigma_2^{-1}\sigma_1^{-1}\sigma_2^{-1}\dots$
($3p$ alternating letters). It represents the braid $\Delta^{-p}$.
By Lemmas 1 and 2, it suffices to prove that the modifications in Figure 3
allow us to transform any
quasipositive graph with the boundary word
$W\Delta_{-p}$ so that every branch point of the new graph is connected by an edge to a point
on $\partial D$.
We say that such branch points are {\it good} and the others are {\it bad}.

We prove this fact by induction on the weight of the graph which we define as
$|R(\Gamma)|$ plus the number of bad branch points. If the weight is zero, then all
branch points are good. Let us prove the required result for a graph $\Gamma$
assuming that it is proven for all graphs of smaller weight. If there are no bad points,
we are done. So, we assume that they exist. Suppose that no modification reduces the weight.
Let $b$ be a bad branch point. Without loss of generality, we may assume that the edge $bv$
incident to $b$ is labeled by 1. In a neighbourhood of this edge, $\Gamma$ is oriented
as in Fig.~1 with $i=1$ and $j=2$ because otherwise modification A in Fig.~3 would reduce 
$|R(\Gamma)|$. Let $P$ be the closure of the component of $D\setminus\Gamma$ that contains $b$.
Let $e_1,\dots,e_n$ be the edges of $\Gamma$ lying on $\partial P$,
numbered in the order of a positive circuit along $\partial P$ starting from $v$.

Let us prove by induction on $i$ that if $e_i$ has a positive (resp. negative)
orientation with respect to $\partial P$, then it is labeled by $2$ (resp. by $1$).
Since this fact contradicts the orientation of $e_n$ in Fig.~1, this will complete
the proof of Theorem 1.
The assertion is true for $i=1$. Suppose that it is true for some $i$. To deduce it for $i+1$,
it suffices to exclude all the cases shown in Fig.~4.

\midinsert
\centerline{  \epsfxsize=77mm\epsfbox{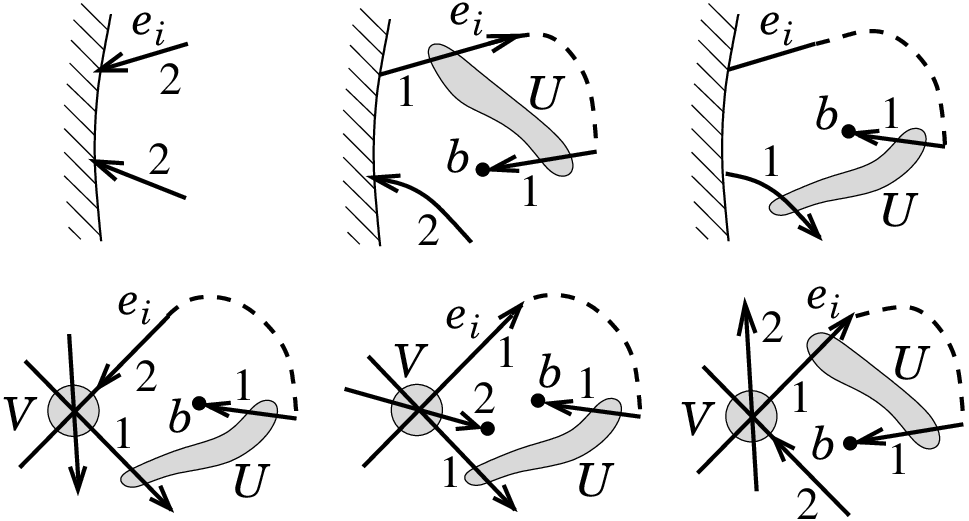} }
\smallskip
\centerline{  {\bf Fig.~4.} }
\endinsert

The left upper case is impossible due to our choice of $\Delta_{-p}$.
The two right upper cases are impossible because modification A
applied to the disk $U$ reduces the number of bad branch points.
The others are impossible because modification A
applied to $U$ followed by modification B applied to $V$ reduces $|R(\Gamma)|$.
Theorem 1 is proven.


\subhead \S3. Proof of Theorem 2
\endsubhead
Let $\Cal A=\{\sigma_0,\sigma_1,\sigma_2\}$ and let $\Cal P$ be the set of
BKL-positive words. Recall that $U=V$ is the equality in $\B_3$ whereas
$U\letterwise V$ is the letterwise coincidence of words.
We follow the convention that if the notation $\letterwise$ is used, then all factors
on both sides of the equality belong to $\Cal P$. We set $\tau:\B_3\to\B_3$ and
$\tau(X)=\delta^{-1}X\delta$. Since $\tau(\Cal A)=\Cal A$, we can define
$\tau:\Cal P\to\Cal P$.

\proclaim{ Lemma  } {\rm(see [1; Theorem 2.7]).}
If $U=V$ for $U,V\in\Cal P$, then $V$ is obtained from $U$
by relations {\rm(2) (}without inserting $\sigma_i^{-1}${\rm)}.
\qed\endproclaim

\proclaim{ Lemma 4 } If $a_1\dots a_n=b_1\dots b_n$,
$a_i,b_i\in\Cal A$, then $(a_1,\dots,a_n)\sim(b_1,\dots,b_n)$.
\qed\endproclaim

\proclaim{ Lemma 5 }
Let $\delta^p=W\equiv AuB$, $u\in\Cal A$. Then there exist $v\in\Cal A$ and $k\ge0$
such that either $A\equiv A_1vA_2$,  $A_2=\delta^k$, $vA_2u=\delta^{k+1}$
or $B\equiv B_1vB_2$,  $B_1=\delta^k$, $uB_1v=\delta^{k+1}$.
\endproclaim

\demo{ Proof } If $p=1$, the statement is obvious.
Assume that it is true for $p-1$. By Lemma , either we have
$W\letterwise(\sigma_2\sigma_1)^p$, or at least one relation can be applied to $W$.
In both cases we have $W\letterwise CxyD$ where $xy=\delta$.
If $A\letterwise C$ or $A\letterwise Cx$, then we set $v=y$ or $v=x$ respectively and we are done.
Otherwise we have either $C\letterwise AuC_1$, or $D\letterwise D_1uB$.
We consider only the former case (the latter one is similar).
Let $E=\tau^{-1}(D)$. Then $AuC_1E=CE=C\delta D\delta^{-1}=W\delta^{-1}=\delta^{p-1}$.
Hence, by the induction hypothesis, we have either
1) $A\letterwise A_1vA_2$, $A_2=\delta^k$, $vA_2u=\delta^{k+1}$, or
2) $C_1\letterwise B_1vC_2$, $B_1=\delta^k$, $uB_1v=\delta^{k+1}$, or
3) $E\letterwise E_1wE_2$, $C_1E_1=\delta^m$, $uC_1E_1w=\delta^{m+1}$.
In Cases 1) and 2), the statement of the lemma is evident. In Case 3), we set
$B_1=C_1xy\tau(E_1)$, $v=\tau(w)$, $B_2=\tau(E_2)$, $k=m+1$.
\qed\enddemo

Let $\Cal I_k(W)=\{I\mid W\setminus I=\delta^p,p\ge0,|I|=k\}$.

\proclaim{ Lemma 6 } Let $W=a_1\dots a_n\in\Cal P$,
$I\in\Cal I_k(W)$. Suppose that $a_i a_{i+1}=\delta$ for some $i$ and
that one of the numbers $i$, $i+1$ belongs to $I$ and the other one does not.
Then there exists $J\in\Cal I_k(W)$ such that
$W_I\sim W_J$ and $\{i,i+1\}\cap J=\varnothing$.
\endproclaim

\demo{ Proof }
Let $i\not\in I$ and $i+1\in I$ (in the case $i\in I$ and $i+1\not\in I$ the proof is similar).
Let $W\letterwise W_1x_1\dots W_kx_kW_{k+1}$ where $x_1,\dots,x_k$ are the letters whose
positions belong to $I$. Then, for some $m$, we have $x_m=a_{i+1}$
and $W_m\letterwise Au$ where $u=a_i$. Let $A_j=W_1\dots W_j$ and $X_j=A_jx_jA_j^{-1}$, $j=1,\dots,k$.
Let $v=a_l$ be the letter in $W\setminus I$ matching $u$ whose existence is proven in
Lemma 5, and let $W_s$ be the subword that contains it.
Let $J=\{l\}\cup(I\setminus\{i+1\})$. We are going to prove that $W_J\sim W_I$.
Let $W_J=(Y_1,\dots,Y_k)$. It is clear that $Y_j=X_j$ when $j<\min(m,s)$.

\smallskip
Case 1. $s\le m$. Then $W_s\dots W_m\letterwise BvDu$, $D=\delta^q$, $vDu=\delta^{q+1}$,
and if, moreover, $s<m$, then $W_s\letterwise BvC$ and $D\letterwise CW_{s+1}\dots W_{m-1}A$.
Since $ux_m=\delta$, it follows that $vDu=Dux_m$, hence
$X_m=A_{s-1}Bv(Du x_m)(vDu)^{-1}B^{-1}A_{s-1}^{-1} =
A_{s-1}BvB^{-1}A_{s-1}^{-1}=Y_s$.
If $j>m$, then we have $Y_j=B_jx_j B_j^{-1}$ where
$B_j=A_{s-1}B(Dux_m)W_{m+1}\dots W_j=
A_{s-1}B(vDu)
W_{m+1}\dots W_j
=A_j$, whence $Y_j=X_j$.
If $s\le j<m$, then
$Y_sY_{j+1}Y_s^{-1}=B_j x_j B_j^{-1}$ where
$B_j= Y_s A_{s-1}BC W_{s+1}\dots W_j =
(A_{s-1}B v B^{-1}A_{s-1}^{-1})A_{s-1}BC
W_{s+1}\dots W_j =
A_{s-1}B v C W_{s+1}\dots W_j$ $= A_j$, whence $Y_sY_{j+1}Y_s^{-1}=X_j$. Thus
$W_J=(X_1,\dots,X_{s-1},X_m,
X_m^{-1}X_sX_m,\,\dots,\,
X_m^{-1}X_{m-1}X_m,
$\par\noindent$
X_{m+1},\dots,X_k)
=\Sigma_s^{-1}\dots\Sigma_{m-1}^{-1}(W_I)$.

\smallskip
Case 2. $s>m$. Then $W_s\letterwise BvC$, $ux_m=\delta$, $D=\delta^q$, $uDv=\delta^{q+1}$ where
$D\letterwise W_{m+1}\dots W_{s-1}B$.  Hence $uDv=ux_mD$ and, by canceling $u$, we obtain
$Dv=x_mD$. Therefore,
$Y_{s-1}=A_m x_m (D v)(x_m D)^{-1}A_m^{-1}=A_m x_m A_m^{-1}=X_m$.
If $j\ge s$, then
$Y_j=B_jx_jB_j^{-1}$ where
$B_j=A_m (x_m D)C W_{s+1}\dots W_j =
A_m (D v)C W_{s+1}\dots W_j = A_j$, whence $Y_j=X_j$.
If $m < j < s$, then
$Y_{j-1}= B_j x_j B_j^{-1}$ where
$B_j=A_m x_m W_{m+1}\dots W_j =
(A_m x_m A_m^{-1}) A_j = X_m A_j$, whence $Y_{j-1}=X_mX_jX_m^{-1}$. Thus
$W_J=
(X_1,\dots,X_{m-1},\,X_m X_{m+1}X_m^{-1},\,\dots,\,X_m X_{s-1}X_m^{-1},\,X_m,X_s,\dots,X_k)=
$\par\noindent$
\Sigma_{s-1}\dots\Sigma_m(W_I)$.
\qed\enddemo

\proclaim{ Lemma 7 } Let $xy=\delta$, $x,y\in\Cal A$, and
 either {\rm(i)} $W\letterwise AxyB$ and $V\letterwise AuvB$ where $uv=\delta$,
 or {\rm(ii)}  $W\letterwise AB$ and $V\letterwise Axy\,\tau(B)$,
 or {\rm(iii)} $W\setminus I\ne 1$, $W\letterwise AxyB$, and $V\letterwise A\,\tau^{-1}(B)$.
Suppose that $I\in\Cal I_k(W)$. Then there exists $J\in\Cal I_k(V)$
such that $V_J\sim W_I$.
\endproclaim

\demo{ Proof } Let $i-1$ be the letter length of $A$ and let $m=|\{i,i+1\}\cap I|$.

(i).
If $m=0$, then $V_I=W_I$.
If $m=2$, then the result follows from Lemma 4. The case when $m=1$ reduces to the case when $m=0$
by Lemma 6.

(ii). $V_J=W_I$ for $J=(I\cap[1,i-1])\cup(2+(I\cap[i,k]))$.

(iii). The case when $m=0$ follows from (ii). The case when $m=1$ reduces to the case when $m=0$
by Lemma 6. Let us consider the case when $m=2$, i.~e., $\{i,i+1\}\subset I$.
Let $W\letterwise a_1\dots a_n$. The condition $W\setminus I\ne1$ means that
either $\{1,\dots,i-1\}\not\subset I$, or $\{i+2,\dots,n\}\not\subset I$.
We consider only the former case (the latter is similar). Let
$l=\max(\{1,\dots,i-1\}\setminus I)$. Without loss of generality we may assume that
$a_l=\sigma_2$. Set $C\letterwise a_1\dots a_{l-1}$, $D\letterwise a_{l+1}\dots,a_{i-1}$,
$E\letterwise\tau(D)$.
Since $Dxy=D\delta=\delta E=\sigma_1\sigma_0 E$, we have
$W_I\sim U_I$ by Lemma 4 where $U\letterwise C\sigma_2\sigma_1\sigma_0 EB$. By
Lemma 6, we have $U_I\sim U_K$ where $\{l,l+1\}\cap K=\varnothing$.
Since $C\sigma_2\tau^{-1}(EB)\letterwise C\sigma_2D\tau^{-1}(B)\letterwise V$,
we have reduced the case when $m=2$ to
the case when $m=1$ with
$(A,xy,B;I)$ being replaced by
$(C\sigma_2,\sigma_1\sigma_0,EB;K)$.
\qed\enddemo

Now Theorem 2 follows from Lemma 7. Indeed, any quasipositive factorization
$(A_jx_jA_j^{-1})_{j=1}^k$ of a given braid can be represented in the form
$W_I$ where $W=W_1x_1\dots W_kx_kW_{k+1}$, $W_j\in\Cal P$,
$A_{j-1}^{-1}A_j=W_j\delta^{-3p_j}$ ($A_0=A_{k+1}=1$) and
the replacements (i)--(iii) of Lemma 7\
allow us to transform $W$ to any given word.


\subhead \S4. Proof of Theorem 3
\endsubhead
Let ${]1,\Delta[}=\{\sigma_1,\,\sigma_2,\,\sigma_1\sigma_2,\,\sigma_2\sigma_1\}$,
$\tau(X)=\Delta^{-1} X\Delta$. For $u\in{]1,\Delta[}$, we denote the first and the last
letter of $u$ by $S(u)$ and $F(u)$ respectively.
Any 3-braid can be written in the form $u_1\dots u_n\Delta^{-p}$
with $u_i\in{]1,\Delta[}$ and $F(u_i)=S(u_{i+1})$ (the {\it right normal form\/}).
Any conjugacy class except
$\sigma_1\Delta^{2m+1}$ and $\sigma_1\sigma_2\Delta^{2m}$ contains an element of this form
such that, moreover, $\tau^p(F(u_n))=S(u_1)$. Let $X$ be such an element.
We identify $u_1,\dots,u_n$
(in this order) with the vertices of a regular polygon $P$.
We define an {\it antisymmetry} of $P$ as a reflection $s$ such that its axis passes though
the midpoints of two sides
$ab$ and $cd$, so that $e(u_i)=2$ for $u_i\in\{a,b,c,d\}$, and $e(u_i)\ne e(s(u_i))$
otherwise.
It is not difficult to derive the following fact from
Theorem 1.

\proclaim{ Lemma 8 } If $X$ is as above and
$e(X)=2$, then the number of orbits of the Hurwitz action
on quasipositive factorizations of $X$ is equal to the number of antisymmetries of $P$.
\qed
\endproclaim

Thus, Theorem 3 follows from:

\proclaim{ Lemma 9 } $P$ has at most two antisymmetries.
\endproclaim

\demo{ Proof } Suppose that $P$ has three different antisymmetries $s_1$,
$s_2$, $s_3$. We number them so that the angle $\alpha$ between the axes of
$s_1$ and $s_2$ is minimal, in particular, $\alpha\le\pi/3$.

\smallskip
Case 1. $\alpha=2\pi/n$. Let $u_{-1}u_0$ and $u_0u_1$ be the invariant sides for $s_1$
and $s_2$. We assume here that the indices are defined mod $n$.
Then the antisymmetricity of $s_1$ and $s_2$ implies
$e(u_0)=e(u_{\pm1})=2$, $e(u_{\pm2})=3-e(u_{\mp1})=1$, $e(u_{\pm3})=3-e(u_{\mp2})=2$, etc.
till the vertices antipodal to $u_{\pm1}$. Since the values of $e$ alternate,
it follows that there is no room for the axis of $s_3$.

\smallskip
Case 2. $\alpha>2\pi/n$.
Let $r=s_2\circ s_1$ (a rotation by $2\alpha$) and let $ab$, $cd$ be the sides invariant by
$s_3$. Let $a^\pm=r^{\pm1}(a)$, $b^\pm=r^{\pm1}(b)$.
The condition $2\pi/n<\alpha\le\pi/3$
implies that
$\{a,b,c,d\}\cap\{a^{\pm},b^{\pm}\}=\varnothing$. Since $s_3(a^+)=b^-$,
it follows that $e(a^+)=1$ or $e(b^-)=1$. We assume that
$e(b^-)=1$ (the case when $e(a^+)=1$ is similar).
Note that if $e(u_i)=1$, then always $e(s_j(u_i))=2$.
Hence $e(s_1(b^-))=2$. Since $s_2(b)=s_2(r(b^-))=s_1(b^-)$ and $e(b)=2$,
it follows that $e(b)=e(s_2(b))=2$. Therefore, $b$ and $s_2(b)$ are consecutive
vertices and the axis of $s_2$ passes between them. Hence the angle between the axes $s_2$ and $s_3$
is equal to $2\pi/n$. This contradicts the minimality of $\alpha$.
\enddemo


\Refs
\def\r{\ref}

\r\no1
\by J.~Birman, K.-H.~Ko, S.-J.~Lee
\paper A new approach to the word and conjugacy problems in the braid groups
\jour  Adv. Math. \vol 139 \yr 1998 \pages 322--353
\endref

\r\no2
\by    S.~Kamada
\paper Surfaces in $\Bbb R\sp 4$ of braid index three are ribbon
\jour  J. of Knot Theory and Ramifications \vol 1 \yr 1992 \pages 137--160
\endref

\r\no3
\by    S.~Yu.~Orevkov
\paper Quasipositivity problem for 3-braids
\jour  Turkish Journal of Math. \vol 28 \yr 2004 \pages 89--93
\endref

\r\no4
\by    S.~Yu.~Orevkov
\paper Algorithmic recognition of quasipositive braids of algebraic length two
\jour J. of Algebra \vol 423 \yr 2015 \pages 1080--1108
\endref

\endRefs
\enddocument